\documentclass[a4paper,12pt]{article}
\usepackage[T2A]{fontenc}
\usepackage[cp1251]{inputenc}
\usepackage{amsmath}
\usepackage{amssymb}
\usepackage[dvips]{graphicx}
\usepackage[english,russian]{babel}

\newtheorem{theorem}{Теорема}[section]
\newtheorem{corollary}[theorem]{Следствие}

\newtheorem{proposition}[theorem]{Предложение}
\newtheorem{definitionhead}[theorem]{Определение}

\newcommand\Proof{\noindent{\bf Доказательство. }}
\newcommand\Endproof{\nopagebreak\strut%
\nopagebreak\hfill\nopagebreak$\Box$\medbreak}

\def\Card{\operatorname{Card}}

\def\Const{\operatorname{Const}}

\begin{document}

\title{{\small УДК 512+519.17+517.987} \\ Описание нормальных базисов
  граничных алгебр и факторных языков медленного роста.}
  
   \author {А.~Я.~Белов, А.~Л.~Чернятьев} \date{}

\maketitle

\abstract{

 Для алгебры $A$ через $V_A(n)$ обозначим размерность
вектрного пространства, порожденного мономами длины не больше $n$.
Пусть $T_A(n)=V_A(n)-V_A(n-1)$. Назовем алгебру
{\em граничной}, если $T_A(n)-n<\Const$.
В работе описываются нормальные базисы для алгебр с
 медленным ростом или граничных алгебр. 

 Пусть $\cal L$~-- факторный язык над конечным алфавитом $\cal A$.  {\it Функция роста} $T_{\cal L}(n)$ есть число подслов длины $n$ в $\cal L$. Мы также описываем факторные языки, такие что $T_{\cal L}(n)\le n+\mbox{const}$.
 \footnote{Работа поддержана грантом РФФИ № 14-01-00548.}}

\medskip
{\bf Ключевые слова:} Нормальный базис, последовательность Штурма, функция роста, мономиальная алгебра, Bergman Gap theorem, факторный язык.

\medskip

\section {Введение.}

Проблематике, связанной с ростом в словах и алгебрах, посвящена обширная литература (см. библиграфию в работах \cite{BBL,L1}).
Для алгебры $A$ через $V_A(n)$ обозначим размерность
вектрного пространства, порожденного мономами длины не больше $n$.
Пусть $T_A(n)=V_A(n)-V_A(n-1)$. Если $A$ однородна, то $T_A(n)$
есть размерность векторного пространства, порожденного мономами
длины ровно $n$.

Известно, что либо $\lim_{n\to
\infty}{(T_A(n)-n)}=-\infty$ (в этом случае  :
{\it либо $\lim{V_A(n)}=C<\infty$ и тогда $\dim A
< \infty$, либо $V_A(n)=O(n)$ и алгебра имеет {\em медленный рост},
либо  $T_A(n)-n < \Const$, либо, наконец, $\lim_{n\to
\infty}{(T_A(n)-n)}=\infty$.} (Bergman Gap Theorem)

В последнем случае для любой функции $\phi(n)\to \infty$ и любой
$\psi(n)=e^{o(n)}$ существует алгебра $A$ такая, что для
бесконечного множества натруальных чисел $n\in L\subset {\mathbb N}\
T_A(n)>\psi(n)$ и для бесконечного множества натуральных чисел $n\in
M\subset {\mathbb N}\ T_A(n)<n+\phi(n)$ (\cite{BBL}) и  рост
может быть хаотичным.

Случай алгебр {\em медленного
роста}, когда $T_A(n)-n < \Const$  исследовался рядом авторов, в частности Дж.~Бергманом и Л.~Смоллом \cite{Be,BBL}.  Назовем алгебру
{\em граничной}, если $T_A(n)-n<\Const$.

Объекты минимального роста представляют значительный интерес.
Для слова $W$ можно определить {\it функцию роста}:
$T_W(n)=\Card F_n(W)$. Ясно, что если $T_W(n)=0$ для какого-то $n$, то $W$ -- конечное
слово.   Можно показать, что если $W$~-- непериодическое слово, то $T_W(1)\ge 2$ и кроме того $T_W(n+1)>T_W(n)$, поэтому $T_W(n)\ge n+1$, если $T_W(n)=n+1$ о получаются знаменитые {\it последовательности Штурма} (см. чрезвычайно красивую классическую теорему \ref{TheorEq}), которым посвящена обширная библиография \cite{BS}. Слова, удовлетворяющие более слабому условию $T_W(n)<n+C$ описаны в теореме \ref{MinGrow}.

Язык $\cal L$ называется {\it факторным} если подслово слова из $\cal L$ снова принадлежит $\cal L$. Язык называется {\it граничным}, если $T_{\cal L}(n)\le n+\mbox{const}$.

Наша цель состоит в
описании нормальных базисов граничных алгебр и граничных факторных языков, что обобщает теорему \ref{MinGrow}. Такое описание сводится к мономиальному случаю \cite{BBL}. Оно дается в теореме \ref{ThMainLast}.

Назовем {\em обструкцией} приводимое (то есть уменьшаемое) слово $u$
такое, что любое его подслово неприводимо. {\it Сверхсловом} в алгебре $A$
(правым, левым, двусторонним) называется сверхслово $W$ такое, что
любое его конечное подслово ненулевое. Аналогично определяется {\em
неприводимое сверхслово} в алгебре $A$. $|w|$ есть длина слова $w$.

Сверхслово $W$ называется {\em рекуррентным}, если каждое его подслово
встречается в нем бесконечно много раз (в случае двустороннего
бесконечного слова, каждое подслово встречается бесконечно много
раз в обоих направлениях). Слово $W$ называется {\em
равномерно-рекуррентным} или ({\em р.р словом}), если оно
рекуррентно и для каждого подслова $v$ существует натуральное
$N(v)$, такое, что для любого подслова $W$ $u$ длины не менее, чем
$N(v)$, $v$ является подсловом $u$.

%Здесь и далее через  $A$ будем обозначать конечный алфавит, то есть
%непустое множество элементов (символов).  Через $A^+$ обозначим множество
%всех конечных последовательностей, символов, или {\em слов}.

%Слово $u$ есть {\em подслово} (или {\em фактор}) слова $w$, если
%существуют слова $p,q\in A^+$ такие, что $w = puq$.
%Если слово $p$ (или $q$) равно $\Lambda$, то $u$ называется {\em
%префиксом} (или {\em суффиксом}) слова $w$.

%Рассмотрим бесконечное (одностороннее или двустороннее) слово $W$ над
%алфавитом $A$. Пусть $v$ -- его подслово и $x\in A$.
%Символ $x$ -- {\it левое (правое) расширение} $v$, если $xv$
%(соотв. $vx$) принадлежит $F(W)$.

%Подслово $v$ называется  {\em левым (правым) специальным
%подсловом}, если для него существуют два или более левых (правых)
%расширения.
%Подслово $v$ называется {\em биспециальным}, если оно является и
%левым, и правым специальнм подсловом одновременно.
%Количество различных левых (правых) расширений подслова
%назовем  {\em левой (правой) валентностью} этого подслова.

%Пусть $W$ -- сверхслово. Для любого подслова $v$ можно
%определить множество {\em возвращаемых} слов $\Ret_W(v)$, а
%именно, слово $u$ -- {\em возвращаемое} для $v$, если $vuv$ --
%подслово $W$ и $v$ -- не является подсловом $u$. Ясно, что для
%рекуррентных слов множество возвращаемых слов $\Ret(v)$ каждого
%подслова $v$ будет непустым, а в случае равномерной рекуррентности
%множество длин слов из $\Ret(v)$ будет ограниченным.

Важным инструментом для описания слова $W$ является {\em графы
подслов}, или {\it графы Рози (Rauzy's graphs)},
которые строятся так: {\em $k$-граф}
слова $W$ -- ориентированный граф, вершины которого
взаимнооднозначно соответствуют подсловам длины $k$ слова $W$, из
вершины $A$ в вершину $B$ ведет стрелка, если в $W$ есть подслово
длины $k+1$, у которого первые $k$ символов -- подслово
соответствующее $A$, последние $k$ символов -- подслово,
соответствующее $B$. таким образом, ребра $k$-графа биективно
соответствуют ($k+1$)-подсловам слова $W$. В этих терминах описываются динамические системы, связанные с перекладываниями отрезков \cite{BeChe}.

%\section{Слова Штурма.}

 Пусть $M$ -- компактное метрическое пространство, $U\subset M$~--- его открытое подмножество, $f:M\to M$ -- гомеоморфизм компакта в себя и $x\in M$ -- начальная точка. По последовательности итераций можно построить бесконечное слово
над бинарным алфавитом: $
w_n=\left\{
\begin{array}{rcl}
   a,\ f^{(n)}(x_0)\in U\\
   b,\ f^{(n)}(x_0)\not\in U\\
\end{array}\right.
$
которое называется {\em эволюцией} точки $x_0$.
Для слов над алфавитом, состоящим из большего числа символов нужно
рассмотреть несколько характеристических множеств: $U_1,\ldots
,U_n$.
Пусть
$\mathbb{S}^1$ -- окружность единичной длины, $U\subset\mathbb{S}^1$
-- дуга длины $\alpha$, $T_\alpha$ -- сдвиг окружности на
иррациональную величину $\alpha$.
Слова, получаемые такими динамическими системами, называют
{\em механическими словами} ({\em mechanical words}).
Имеется классическая

\begin{theorem}[Теорема эквивалентности \cite{BS}]\label{TheorEq}
Следующие условия на свер\-х\-с\-ло\-во $W$ почти эквивалентны
(т.е. точностью до счетного множества, исключительные последовательности также описаны):

\begin{enumerate}
\item Слово $W$ является {\it словом Штурма}, т.е. имеет функцию сложности $T_W(n)=n+1$).

\item Слово $W$ {\em сбалансированно} (т.е.  для любых его двух
подслов $u,v$ одинаковой длины выполняется неравенство:
$
||u|_a-|v|_a| \leq 1
$).

\item Слово $W$ является {\it механическим}.%, то есть порождается системой
  %$(\mathbb{S}^1,U,T_\alpha)$ с иррациональным $\alpha$.
\end{enumerate}
\end{theorem}

Естественным обобщением слов Штурма являются  слова с минимальной функцией роста, то есть с функцией роста,
удовлетворяющей соотношению $F_W(n+1)-F_W(n)=1$ при всех
достаточно больших $n$.
В работе \cite{Che1} доказана следующая

\begin{theorem}\label{MinGrow}
Пусть $W$ -- рекуррентное слово над произвольным конечным
алфавитом $A$. Тогда следующие условия на слово $W$ эквивалентны:
\begin{enumerate}

\item Существует такое натуральное $N$, что функция сложности
слова $W$ равна $T_W(n)=n+K$, для $n\geq N$  и некоторого
постоянного натурального $K$.

\item Существуют такое иррациональное $\alpha$ и целые $n_1,n_2,
\ldots,n_m$, что слово $W$ порождается динамической системой
$$
(\mathbb{S}^1,T_\alpha, I_{a_1},I_{a_2},\ldots ,I_{a_n},x),
$$
где
$T\alpha$ -- сдвиг окружности на иррациональную величину $\alpha$,
$I_{a_i}$ -- объединение дуг вида $(n_j\alpha,n_{j+1}\alpha)$.

\end{enumerate}

\end{theorem}

\section{Основные теоремы.}

\begin{theorem}[О нормальных базисах]\label{ThMainLast}
Пусть $A$~-- граничная алгебра  с образующими  $a_1,\ldots,a_s$.
Тогда имеют место один из двух случаев.

{\bf Случай 1.} Алгебра $A$ не содержит равномерно-рекуррентного
непериодического сверхслова. В этом случае нормальный базис алгебры
$A$ состоит из множества подслов следующего множества слов:

\begin{enumerate}
\item Одно слово вида $W=u^{\infty/2}cv^{\infty/2} \neq u^{\infty}$

\item Произвольный конечный набор $\mu$ конечных слов

\item Множество слов вида $u_i^{\infty/2}c_i$,\ $i=1,\ldots, r_1$

\item Множество слов вида $d_iv_i^{\infty/2}$,\ $i=1,\ldots, r_2$

\item Множество слов вида $e_j{(R_j)}^{k} f_j$,\ $k\in \mathbb{K}_j \subseteq
\mathbb{N}$, $j=1,\ldots,r_3$

\item Множество слов вида $W_{\alpha} = E_{\alpha}u^{n_\alpha}cv^{m_\alpha}
F_\alpha$. При этом существуют такое $c>0$, что для любого $k$
количество слов $W_\alpha$ длины $k$ меньше $c$.
\end{enumerate}

{\bf Случай 2.} Алгебра $A$ содержит равномерно-рекуррентное
непериодическое сверхслово $W$. В этом случае нормальный базис
алгебры $A$ состоит из множества подслов следующего семейства слов,
включающее в себя:

\begin{enumerate}
\item Некоторое равномерно рекуррентное слово $W$ с функцией роста
$T_W(n)=n+const$ для всех достаточно больших $n$. Описание таких
слов дано в теореме \ref{MinGrow}.

\item Произвольный конечный набор $\mu$ конечных слов

\item Множество слов вида $u_i^{\infty/2}c_i$, $i=1,\ldots,s_1$

\item Множество слов вида $d_iv_i^{\infty/2}$, $i=1,\ldots,s_2$

\item Множество слов вида $e_j(R_j)^{k}f_j$, $k\in \mathbb{K}_j \subseteq
\mathbb{N}$,\ $j=1,\ldots,s_3$

\item Множество слов вида $L_iO_iW_i$,\ $i=1,\ldots,k_1$. При этом
$W_i$ -- сверхслово, эквивалентное $W$ и $O_iW_i$ имеют вхождение
только одной обструкции (а, именно, $O_i$).

\item Множество слов вида $W'_jO'_jL'_j$, $j=1,\ldots,k_2$. При этом
$W'_j$ -- сверхслово, эквивалентное $W$ и $W'_jO'_j$ имеют вхождение
только одной обструкции (а, именно, $O'_j$).

\item Конечное множество серий вида: ${h^1}_iT_i{h^2}_i$,
$i=1,\ldots,s$. При этом:

a) слово $T_i$ содержит вхождение ровно двух обструкций
$O_{i_1},O_{i_2}$ (возможно, перекрывающихся)

b) Для некоторого $c>0$ $|h^1_i|+|h^2_i|<c$ при всех $i$.

с) Существует $m>0$ такое, что для любого $k$ имеется не более $k$
подслов длины $m$, вида $(8)$ и не являющихся подсловами слова $W$.

\end{enumerate}
\end{theorem}

\Proof Алгебру $A$ считаем мономиальной. Назовем слово $v$ алгебры
$A$ {\em хорошим}, если для любого $n$ существуют сколь угодно
длинные слова $w_1,w_2$, $|w_1|>n,|w_2|>n$ такие, что $w_1vw_2$ есть
подслово алгебры $A$. Обозначим $T_{RL}(n)$  количество хороших слов
длины $n$. Известно, что если $T_{RL}(n)=T_{RL}(n+1)$ при некотором
$n$, то алгебра $A$ имеет медленный рост (см. \cite{BBL}).  В силу
граничности алгебры $A$ $T_{RL}(n)\geq T_{RL}(n+1)+1$,  при всех
достаточно больших $n$ неравенство превращается в равенство (иначе
$\lim_{n\to \infty} (T_{RL}(n)-n)= \infty$ и, как следствие
$\lim_{n\to \infty} (T_{A}(n)-n)= \infty$).

При этом граф Рози имеет развилку и, как следствие, два цикла,
эволюция графа Рози устроена следующим образом: либо граф теряет
сильную связность и имеет вид a). В этом случае его дальнейшая
эволюция однозначна, она отвечает слову вида $u^{\infty/2} c
v^{\infty/2}$ и мы имеем случай $1$, либо граф Рози все время
остается сильно связным. Тогда эволюция связной компоненты, в
которой есть развилка, асимптотически эквивалентна эволюции графа
Рози некоторого слова Штурма. Если оно имеет вид  $u^{\infty/2} c
u^{\infty/2}$, то имеет место случай $1$, иначе оно
равномерно-рекуррентно и имеет место случай $2$.

В случае 1 все обструкции для слова $u^{\infty/2} c u^{\infty/2}$
имеют ограниченную длину (см. \cite{BBL}). Можно сделать следующее

\medskip
{\bf Наблюдение.}\ {\em Количество ненулевых слов длины $n$, не
являющихся подсловами слова $u^{\infty/2}cu^{\infty/2}$ не
превосходит константы (не зависящей от $n$).}
\medskip

\begin{proposition}[\cite{BBL}]
Пусть $SW=WT$. Тогда $W$ имеет вид: $s^k s_1$, где $s_1$ -- начало
слова $s$.
\end{proposition}

Из данного предложения и только что сделанного наблюдения вытекает

\begin{corollary}
Пусть $|u|=k$, $|v_1|=|v_2|=l$. Тогда либо количество подслов длины
$k+m$, (где $m\le k$) слова $v_1uv_2$ не менее $m+1$, либо
$u=s^ks'$, для некоторого слова $s$, при этом $s'$ -- начало $s$ и
$v_1=v'_1s'$.
\end{corollary}

Из данного следствия и сделаного наблюдения получается

\begin{proposition}
Существует константа $K$, зависящая только от граничной алгебры $A$,
такая, что для любой обструкции $O$ в слове $W$ либо при некотором
$m$ для любых $v_1,v_2$, $|v_1|\geq m$, $|v_2|\geq m$ $v_1Ov_2$
является нулевым словом алгебры (число $m$ не зависит от выбора
обструкции), либо $|v|\leq K$.
\end{proposition}

Из данного предложения следует, что словами алгебры $A$ являются
либо слова, содержащие не более двух обструкций, причем каждая
обструкция находится на ограниченном расстоянии от одного из концов,
либо подслова слов вида $R_iu^k_iT_i$. А все такие типы слов описаны
в условии теоремы \ref{ThMainLast}. \Endproof

Мы также доказали следующее утверждение 

\begin{theorem}[О граничных факторных языках]\label{ThMainLast1}
Любой граничный факторный язык есть язык нормальных слов граничной алгебры. Множество его слов описывается в предыдущей теореме. 
\end{theorem}


\begin{thebibliography}{MH1}

\bibitem{BBL} А.Я.Белов, В.В.Борисенко, В.Н.Латышев, Мономиальные
алгебры // Итоги науки и техники. Совр. Мат. Прил. Тем. Обзоры т.
26 (алг. 4), М. 2002. 35-214.

%\bibitem{BK} А.Белов, Г.Кондаков, Обратные задачи символической
%динамики // Фундаментальная и прикладная математика, Т1, №1, 1995

\bibitem{BeChe} Белов А. Я. Чернятьев А.Л., Слова медленного роста и
перекладывания отрезков // Успехи Мат. Наук, 2008, 63:1(379), стр.
159--160.

\bibitem{Che1}
Чернятьев А.~Л. Слова с минимальной функцией роста. // Вестник МГУб Матем. и Мех,
IAO, том 6, стр. 42--44, 2008 г.

%\item{Che} Чернятьев А. Л.,Сбалансированные слова и символическая
%динамика, Фундаментальная и прикладная математика, Том 13, выпуск 5,
%2007 г., стр. 213--224.

%\bibitem{AR} P. Arnoux and G. Rauzy [1991],
%Representation geometrique des suites the complexite $2n + 1$ //
%Bull. Soc. Math. France 119, 199-215.

\bibitem{Be} G.M. Bergman, A note on growth functions of
algebras and semigroups // Research Note, University of
California, Berkeley, 1978.

\bibitem{BS} J. Berstel, P. S\'{e}\'{e}bold, Sturmian words, in: M. Lothaire (Ed.) // Algebraic Combinatorics on Words,
Encyclopedia of Mathematics and Its Applications, Vol. 90,
Cambridge University Press, Cambridge, 2002 (Chap. 2).

\bibitem{L1} M. Lothaire, Combinatorics on Words // Encyclopedia of
Mathematics and its Applications, Addison-Wesley, Reading, MA,
1983, Vol. 17.

%\bibitem{Ra2} G. Rauzy, Suites \'a termes dans un alphabet fini //
%In {\it S\'emin. Th\'eorie des Nombres,} p. 25-01-25-16,
%1982-1983, Expos\'e No 25.

%\bibitem{V} L. Vuillon, Balanced words // Bull. Belg. Math. Soc. Simon Stevin 10 (2003), no. 5, 787-805


\end{thebibliography}
\end{document}